\documentclass{article}
\usepackage[utf8]{inputenc}
\usepackage{amsmath}
\title{A generalization of Iseki Formula and the transformation of $\theta_1(z,\tau)$}
\author{Maher Me'meh $\&$ Ali Saraeb }
\date{}

\begin{document}

\maketitle
\textbf{Abstract} In this paper we give a generalization of Iseki's formula and use it to prove the transformation law of $\theta_1(z,\tau)$.
\vspace{0.3cm}
\\
\begin{center}
INTRODUCTION
\end{center}
The Dedekind Eta function, defined as
\[\eta(\tau)=e^{\pi i \tau/12}\prod_{n=1}^{\infty}(1-e^{2\pi i n \tau}),\]
plays an important role in the study of Modular and Jacobi forms. Its transformation over a matrix $A\in \Gamma$, where  $\Gamma$ is the full modular group, is given by

\[\eta(A\tau)=\epsilon(A)(-i(c\tau+d))^{1/2}\eta(\tau).\]
One obtains the Dedekind eta-character ``$\epsilon(A)"$ which is defined as, (check, \cite{Rademacher},III)
\[\displaystyle \epsilon(A)=\left\lbrace \begin{array}{ll}
\left(\frac{d}{c}\right)i^{(1-c)/2}e^{(\pi i/12)(bd(1-c^2)+c(a+d))}~~~~\text{if}~ c~odd \\
\\
\left(\frac{c}{d}\right)e^{\pi d i/4}i^{(1-d)/2} e^{(\pi i /12)(ac(1-d^2)+d(b-c))}~~~~\text{if}~d~odd,
\end{array}\right.\]
where $(d,c)=1$.\\
However, it turns out that computations done using this definition of the eta-character can get really messy. On the other side, Sho Iseki proved in 1952 (\cite{apostol}, III) the transformation law of the eta function using a functional equation which will be addressed below. Using his proof he was able to write the eta-character using Dedekind sums which proved to be much easier in terms of computations. The eta-character turns out to be

\[\epsilon(A)=exp\left(\pi i \left(\frac{a+d}{12c}-s(d,c)\right)\right),\]
where \[ \displaystyle s(d,c)= \displaystyle\sum_{r=1}^{c-1}\frac{r}{c}\left(\frac{dr}{c}\left[\frac{dr}{c}\right]-\frac{1}{2}\right)\] is the Dedekind sum and $(d,c)=1$.\\
Since the eta-character appears significantly in the transformation laws of Jacobi theta functions and Jacobi forms (\cite{Rademacher},X), we generalize Iseki's proof of the eta function and apply the generalization to Jacobi theta function $\theta_1(z,\tau)$, which is defined as 
\begin{equation}
\theta_1(z,\tau) = -i w q^{1/4} \prod_{n=1}^{\infty} (1-q^{2n})(1- w^2 q^{2n})(1- w^{-2} q^{2n-2}),
\end{equation}
where $w=e^{\pi i z}$ and $q=e^{\pi i \tau}$, $z\in \mathbf{C}~and~\tau\in \mathbf{H}$. We first generalize Iseki's functional equation,
\[\Lambda(\alpha,\beta,z)=\Lambda(1-\beta,\alpha, z^{-1})+g_0(\alpha, \beta,z)\]
to four variables using methods from Fourier analysis then we employ this tool to prove the transformation law of $\theta_1$.

\begin{center}
GENERALIZATION OF ISEKI'S FORMULA
\end{center}
\textbf{Theorem 1.} If $Re(w)>0$,  \(0 <  \alpha < 1\) , $\theta$ is real, and $0<\beta +\theta<1$ , then
\begin{align}
\Lambda(\alpha,\beta,w,\theta)=\Lambda(1-\beta,\alpha,w^{-1},-i\theta/w)+g_0(\alpha,\beta,w,\theta),
\end{align}
where \\ 
\begin{align}
~~~~~~~\Lambda (\alpha,\beta,w,\theta) =-\sum_{n=0}^{\infty} log(1-e^{2\pi i \theta}e^{-2\pi((n+\alpha)w-i\beta)})+log(1-e^{-2\pi i \theta}e^{-2\pi((n+1-\alpha)w+i\beta)}),\end{align}
with
\begin{align}
g_0(\alpha,\beta,w,\theta)=\frac{\pi}{w}B_2(\beta+\theta)-\pi w B_2(\alpha)+2\pi i B_1(\alpha)B_1(\beta+\theta),
\end{align} 
and $B_n$ is the $n^{th}$ Bernoulli polynomial. \\ \\
\textit{Proof.} \\ \\
The proof utilizes the following well-known identities from Fourier analysis: \\
\begin{align} \displaystyle
\frac{e^{2\pi  m \alpha w}}{ 1-e^{2\pi m w}} +\frac{1}{ 2 \pi  w m}= \frac{1}{2 \pi i}\sum_{n=-\infty }^{\infty } \frac{e^{2\pi i \alpha n}}{wmi +n}. \\ \nonumber \end{align} Replacing $w$ by $w^{-1}$ and $m$ and $-m$, we get\\ 
\begin{align}
\frac{e^{-2\pi  m \alpha w^{-1}}}{ 1-e^{-2\pi m w^{-1}}} -\frac{1}{ 2 \pi  w^{-1} m}= -\frac{1}{2 \pi i}\sum_{n=-\infty }^{\infty } \frac{e^{2\pi i \alpha n}}{w^{-1}mi -n}. \\ \nonumber \end{align}
We also have
\begin{align}
\frac{1}{m(wmi-n)}= -\frac{1}{mn} +\frac{w}{n(ni+wm)}.  
\end{align} 
We first observe that (3) can be rewritten as follows \\ 
\begin{align}
~~~~~~~\Lambda (\alpha,\beta,w,\theta) =&  \sum_{m=1}^{\infty} \frac{e^{2\pi i m \beta }}{m}  \frac{e^{-2\pi m \alpha w}}{1- e^{-2 \pi m w}} e^{2\pi i m \theta}  - \sum_{m=1}^{\infty} \frac{e^{-2\pi i m \beta }}{m}  \frac{e^{2\pi m \alpha w}}{1- e^{2 \pi m w}} e^{-2\pi i m \theta} \nonumber\\ =& - \sum_{\substack{m=- \infty  \\ m \neq 0}  }^{\infty} \frac{e^{-2\pi i m \beta }}{m}  \frac{e^{2\pi m \alpha w}}{1- e^{2 \pi m w}} e^{-2\pi i m \theta} .
\end{align}
Multiplying both sides of equation (5) by $\frac{-1}{2\pi mi} e^{-2\pi i m  \beta } e^{-2\pi i m \theta }$ and then summing from $m=-\infty$ to $+\infty$, we rewrite (8) as follows
\begin{align}\displaystyle
~~~~~~~\Lambda (\alpha,\beta,w,\theta) =& -\frac{1}{2\pi i} \sum_{m=-\infty }^{\infty } \sum_{n=-\infty }^{\infty } \frac{e^{-2\pi i m \beta }}{m}   \frac{e^{2\pi i n \alpha }}{wmi +n} e^{-2\pi i m  \theta } + \frac{1}{2 \pi w} \sum_{m=-\infty }^{\infty } \frac{e^{-2\pi i m (\beta+ \theta) } }{m^2} \nonumber\\ 
=& -\frac{1}{2\pi i} \sum_{m=-\infty }^{\infty } \sum_{n=-\infty }^{\infty } \frac{e^{-2\pi i m \beta  }}{m}   \frac{e^{-2\pi i n\alpha }}{wmi-n}e^{-2\pi i m  \theta } + \frac{1}{2 \pi w} \sum_{m=-\infty }^{\infty } \frac{e^{2\pi i m (\beta+ \theta) } }{m^2} \nonumber \\ =& -\frac{1}{2\pi i}
\sum_{m=-\infty }^{\infty } (A_{m,n}(\alpha, \beta , w, \theta)) + \frac{1}{2 \pi w} F_2(\beta+\theta )  ,
\end{align}
where 
\begin{align} \displaystyle 
A_{m,n}(\alpha, \beta , w, \theta)= \sum_{n=-\infty }^{\infty } \frac{e^{-2\pi i m \beta  }}{m}   \frac{e^{-2\pi i n\alpha }}{wmi-n}e^{-2\pi i m  \theta } \\ \nonumber
\end{align}
and
\begin{align} \displaystyle 
F_n(x)= \sum_{m=-\infty }^{\infty } \frac{e^{2\pi i m x } }{m^n} = \begin{cases} \displaystyle \sum_{m=-\infty }^{\infty } \frac{e^{-2\pi i m x } }{m^n} \text{ if n is even. } \\ \\
\displaystyle  -\sum_{m=-\infty }^{\infty } \frac{e^{-2\pi i m x } }{m^n} \text{ if n is odd. }
\end{cases} \end{align}\\  \\
\( \displaystyle
\text{Using (7), (10) becomes as follows} \\ \\
~~~~~~~A_{m,n}(\alpha, \beta , w, \theta)=   -\frac{e^{-2\pi i m (\beta+ \theta) }}{m} \sum_{n=-\infty }^{\infty } \frac{e^{-2\pi i n \alpha  }}{n}    + \sum_{n=-\infty }^{\infty }   e^{-2\pi i m \beta  } \frac{ e^{-2\pi i n \alpha }}{n(\frac{ni}{w}+m)}e^{-2\pi i m \theta  }. \\ \\ \) \\
By using $F_1(x)= -2 \pi i B_1(x)$, $F_2(x)= \frac{-(2 \pi i)^2}{2!} B_2(x) $ , (7), (10), and (11) and by carefully manipulating the signs of $m$ and $n$ in the summands,  we observe that   \\ \\
\(\displaystyle
~~~~~~~-\frac{1}{2 \pi i}\sum_{m=-\infty }^{\infty } A_{m,n}(\alpha, \beta , w, \theta) \)
\begin{align}
~~~~~~~~~~~~~~=&  \frac{1}{2 \pi i}\sum_{m=-\infty }^{\infty } \frac{e^{-2\pi i m (\beta+ \theta) }}{m} \sum_{n=-\infty }^{\infty } \frac{e^{-2\pi i n \alpha }}{n} \nonumber- \frac{1}{2 \pi i}\sum_{m=-\infty }^{\infty } \sum_{n=-\infty }^{\infty }   e^{-2\pi i m (\beta+ \theta) } \frac{ e^{-2\pi i \alpha n}}{n(\frac{ni}{w}+m)} \nonumber\\ \nonumber\\ 
=&  \frac{1}{2 \pi i}F_1(\beta+ \theta) F_1(\alpha) - \frac{1}{2 \pi i} \sum_{m=-\infty }^{\infty } \sum_{n=-\infty }^{\infty }   e^{-2\pi i n (\beta+ \theta) } \frac{ e^{-2\pi i m\alpha }}{m(\frac{mi}{z}+n)} \nonumber\\ \nonumber\\ 
=&  2 \pi i B_1(\beta+ \theta) B_1(\alpha) -\frac{1}{2 \pi i}\sum_{m=-\infty }^{\infty } \sum_{n=-\infty }^{\infty }   e^{2\pi i n (\beta+ \theta) } \frac{ e^{-2\pi i m \alpha }}{m(\frac{mi}{w}-n)} \nonumber\\ \nonumber\\
=& 2 \pi i B_1(\beta+ \theta) B_1(\alpha) + \sum_{m=-\infty }^{\infty } \frac{e^{-2\pi  m (\beta+ \theta) w^{-1}}}{ 1-e^{-2\pi m w^{-1}}}\frac{ e^{-2\pi i m \alpha }}{m} \nonumber-\frac{w}{2 \pi } \sum_{m=-\infty }^{\infty }\frac{ e^{-2\pi i m \alpha }}{m^2}  \nonumber\\ 
=& 2 \pi i B_1(\beta+ \theta) B_1(\alpha) -  \sum_{m=-\infty }^{\infty } \frac{ e^{-2\pi i m \alpha }}{m}\frac{e^{-2\pi  m \beta w^{-1}}}{ 1-e^{2\pi m w^{-1}}} e^{-2\pi  m  \theta w^{-1}} \nonumber -\frac{w}{2 \pi } \sum_{m=-\infty }^{\infty }\frac{ e^{-2\pi i m \alpha }}{m^2}  \nonumber\\
=& 2 \pi i B_1(\beta+ \theta) B_1(\alpha) 
+ \Lambda(1-\beta,\alpha,w^{-1},-i\theta/w) -\pi w B_2(\alpha). \\\nonumber
\end{align}
Plugging (12) in (9), we get that \\ 
\begin{align*}
~~~~~~~\Lambda (\alpha,\beta,w,\theta) = \Lambda(1-\beta,\alpha,w^{-1},-i\theta/w)  +\frac{\pi}{w}B_2(\beta + \theta) -\pi w B_2(\alpha)+ 2 \pi i B_1(\beta+ \theta) B_1(\alpha) .
\end{align*}
This completes the proof of theorem 1. \\


We now use Theorem 1. to prove the transformation law for $\theta_1$ under the elements of the full modular group $\Gamma$.
\vspace{0.3cm}
\\
$\textbf{Theorem 2.}$ For $\tau\in \mathbf{H}$ and $z\in \mathbf{C}$ we have
\begin{equation}
\theta_1\left(\frac{z}{c\tau+d},\frac{a\tau+b}{c\tau+d}\right)=\epsilon_1(A)\left(-i(c\tau+d)\right)^{1/2}e^{\frac{\pi i c z^2}{c\tau+d}}\theta_1(z,\tau).
\end{equation}
Here $\epsilon$ appears in the transformation law of the Dedekind eta function as mentioned in the introduction, where again
\[\epsilon(A)=exp(\pi i (\frac{a+d}{12c}+s(-d,c))),\]
and
\[\displaystyle s(h,k)=\sum_{r=1}^{k-1}\frac{r}{k}\left(\frac{hr}{k}-\left[\frac{hr}{k}\right]-\frac{1}{2}\right)\] is the Dedekind sum for $k>0$ and $(k,h)=1$. \\
Hence \[\epsilon_1(A)=-i\epsilon^3=-i.exp\left(3\pi i\left(\frac{a+d}{12c}+s(-d,c)\right)\right).\]
Having (13) is equivalent to proving 
\begin{align*}
~~~~~~~log\left(\theta_1\left(\frac{z}{c\tau+d},\frac{a\tau+b}{c\tau+d}\right)\right)=log(\epsilon_1(A))+\frac{1}{2}log(-i(c\tau+d))+\frac{\pi i c z^2}{c\tau+d}+log\left(\theta_1\left(z,\tau\right)\right).
\end{align*}
Note that using the definition of Dedekind eta function, we have \[\displaystyle \prod_{n=1}^{\infty}(1-q^{2n})=\eta(\tau).e^{\frac{-\pi i \tau}{12}}.\]
Hence (1) becomes
\begin{equation}
\theta_1(z,\tau) = -i w q^{1/4}\left(\eta(\tau)e^{{\frac{-\pi i \tau}{12}}}\right) \prod_{n=1}^{\infty}(1- w^2 q^{2n})(1- w^{-2} q^{2n-2}).
\end{equation}
So from (14) and (15), we have
\vspace{0.25cm}
\\
$\displaystyle 
~~~~~~~~log(-ie^{\pi i z}e^{\frac{\pi i}{4}(\frac{a\tau+b}{c\tau+d}})+log(\eta(\frac{a\tau+b}{c\tau+d}))-\frac{\pi i}{12}\left(\frac{a\tau+b}{c\tau+d}\right)\\~~~~~+\sum_{n=1}^{\infty}log(1-e^{\frac{2\pi i z}{c\tau+d}}e^{2n\pi i (\frac{a\tau+b}{c\tau+d})}) +\sum_{n=1}^{\infty}log(1-e^{\frac{-2\pi iz}{c\tau+d}}e^{2(n-1)\pi i (\frac{a\tau+b}{c\tau+d})})\\~~~~~=log\left(-i.exp\left(3\pi i(\frac{a+d}{12c}+s(-d,c))\right)\right) +\frac{1}{2}log(-i(c\tau+d))+\frac{\pi i c z^2}{c\tau+d}\\~~~~~ + log(-ie^{\pi i z} e^{\frac{\pi i \tau}{4}})+log(\eta(\tau))-\frac{\pi i \tau}{12}$ \\$\displaystyle ~~~~~+\sum_{n=1}^{\infty} log(1-e^{2\pi i z}e^{2n\pi i \tau})+\sum_{n=1}^{\infty}log(1-e^{-2\pi i z}e^{2(n-1)\pi i \tau}).$
\vspace{0.3cm}
\\
Using the fact that, (check (\cite{apostol},III))
\[~~~~~~~~log\left(\eta\left(\frac{a\tau+b}{c\tau+d}\right)\right)=\pi i \left(\frac{a+d}{12c}\right)+\pi i s(-d,c)+\frac{1}{2}log(-i (c\tau+d))+log(\eta(\tau)),\] we obtain
\vspace{0.3cm}
\\
$\displaystyle ~~~~~~\frac{\pi iz}{c\tau+d} +~\frac{\pi i}{6}\left(\frac{a\tau+b}{c\tau+d}\right)+\sum_{n=1}^{\infty}log(1-e^{\frac{2\pi i z}{c\tau+d}}e^{2n\pi i (\frac{a\tau+b}{c\tau+d})}) +\sum_{n=1}^{\infty}log(1-e^{\frac{-2\pi iz}{c\tau+d}}e^{2(n-1)\pi i (\frac{a\tau+b}{c\tau+d})})\\~~~~~=log(-i)+2\pi i(\frac{a+d}{12c})+2\pi is(-d,c))+\frac{\pi i c z^2}{c\tau+d} + \pi i z +\frac{\pi i \tau}{6}\\ ~~~~~+\sum_{n=1}^{\infty} log(1-e^{2\pi i z}e^{2n\pi i \tau})+\sum_{n=1}^{\infty}log(1-e^{-2\pi i z}e^{2(n-1)\pi i \tau}).$
\vspace{0.3cm}
\\
Relocating the terms we get
\vspace{0.25cm}
\\
$\displaystyle ~~~~~~\sum_{n=1}^{\infty}log(1-e^{\frac{2\pi i z}{c\tau+d}}e^{2n\pi i (\frac{a\tau+b}{c\tau+d})}) +\sum_{n=1}^{\infty}log(1-e^{\frac{-2\pi iz}{c\tau+d}}e^{2(n-1)\pi i (\frac{a\tau+b}{c\tau+d})}) \\ ~~~~~~~~=\sum_{n=1}^{\infty} log(1-e^{2\pi i z}e^{2n\pi i \tau})+\sum_{n=1}^{\infty}log(1-e^{-2\pi i z}e^{2(n-1)\pi i \tau})+~log(-i)+~2\pi i(\frac{a+d}{12c})+~2\pi is(-d,c) \\ ~~~~~~+\frac{\pi i c z^2}{c\tau+d} + \frac{\pi iz}{c\tau+d}+ \pi i z +\frac{\pi i}{6}\left(\tau-\frac{a\tau+b}{c\tau+d}\right).~~~~~~~~~~~~~~~~~~~~~~~~~(16)$
\vspace{0.4cm}
\\
Now we introduce a classical change of variable, we set
\[-i(c\tau+d)=v~~~~~~~~a=H, c=k~and~h=-d,\]
such that $Hh\equiv -1~(mod~k).$
\\
Using this change of variable, we have
\[\tau=\frac{iv+h}{k}~~~and~~~~\frac{a\tau+b}{c\tau+d}=\frac{1}{k}\left(H+\frac{i}{v}\right).\]
Moreover,
\[\frac{\pi i}{6}\left(\tau-\frac{a\tau+b}{c\tau+d}\right)=-2\pi i \left(\frac{a+d}{12c}\right)-\frac{\pi}{6k}\left(v-\frac{1}{v}\right).\]
\vspace{0.3cm}
\\
Plugging in (16), where $log(-i)=-\frac{\pi i}{2}$, we obtain the desired functional equation\\
$\displaystyle ~~~~~~~~~~~\sum_{n=1}^{\infty} log(1-e^{\frac{2 \pi z}{v}}e^{\frac{2n\pi i}{k}(H+\frac{i}{v})}) +log(1-e^{-\frac{2 \pi z}{v}}e^{\frac{2(n-1)\pi i}{k}(H+\frac{i}{v})}) \\ \\ ~~~~~~~~~~~= \sum_{n=1}^{\infty}  log(1-e^{2 \pi i z}e^{\frac{2n\pi i}{k}(h+iv)})  +log(1-e^{-2 \pi i z}e^{\frac{2(n-1)\pi i}{k}(h+iv)})  \\ \\
~~~~~~~~~~~+ 2\pi i s(h,k) -\frac{\pi }{6k} (v-\frac{1}{v}) -\frac{\pi i }{2} + \frac{\pi k z^2 }{v} + \pi i z - \frac{\pi z}{v}. ~~~~~~~~~~~(17)$
\vspace{0.25cm}
\\
We now follow Iseki's proof closely and we let \[ \beta=\frac{\phi}{k}~~ where~ ~1\leq \phi\leq k-1.\]\\
However from theorem (1), $0<\beta+\theta<1$ for which one can easily prove that it is equivalent to having $\displaystyle 0<\theta<\frac{1}{k}$.\\
We prove first the case when $k=1$ and then for any integer $k>0$. We then extend the result to the whole plane using analytic continuation except at the endpoints where we treat them separately.\\
Using theorem (1), we have

\[\Lambda(\alpha,\beta,w,\theta)=\Lambda(1-\beta,\alpha,1/w,-i\theta/w)+g_0(\alpha,\beta,w,\theta).~~~~~~~~~~~~(18)\]
For $k=1$, we have $\beta=0$, then we let $\alpha\to 1$ to obtain from (17)\\
$\displaystyle ~~~~~~~~~~\sum_{n=1}^{\infty} log(1-e^{\frac{2 \pi z}{v}}e^{2n\pi i(H+\frac{i}{v})}) +log(1-e^{-\frac{2 \pi z}{v}}e^{2(n-1)\pi i(H+\frac{i}{v})})$\\
$\displaystyle ~~~~~~~~~~=\sum_{n=1}^{\infty}  log(1-e^{2 \pi i z}e^{2n\pi i(h+iv)})  +log(1-e^{-2 \pi i z}e^{2(n-1)\pi i(h+iv)})\\~~~~~~~~~~~~~~~~~~ -\frac{\pi }{6} (v-\frac{1}{v})-\frac{\pi i }{2} + \frac{\pi z^2 }{v}- \frac{\pi z}{v} + \pi i z $.\
\vspace{0.3cm}
\\
Note that $e^{2n\pi i H}=1$, so we end up with 
\vspace{0.3cm}
\\
$\displaystyle ~~~~~~~~\sum_{n=1}^{\infty} log(1-e^{\frac{2 \pi z}{v}}e^{\frac{-2n\pi}{v})}) +log(1-e^{-\frac{2 \pi z}{v}}e^{\frac{-2(n-1)\pi }{v})}) \displaystyle =\sum_{n=1}^{\infty}  log(1-e^{2 \pi i z}e^{-2n\pi v)})  \\~~~~~~+\sum_{n=1}^{\infty}log(1-e^{-2 \pi i z}e^{-2(n-1)\pi v)})\displaystyle  -\frac{\pi }{6} \left(v-\frac{1}{v}\right) -\frac{\pi i }{2} + \frac{\pi z^2 }{v} + \pi i z -\frac{\pi z}{v}.~~~~~~~~~~~~~~~~~(19)$
\vspace{0.25cm}
\\
Using the fact that $\beta=0$ and $\alpha\to 1$, from (18), we see that\\
$\displaystyle ~~~~~-\sum_{n=0}^{\infty} log(1- e^{-2 \pi i \theta } e^{-2\pi (n)w )})+log(1- e^{2 \pi i \theta } e^{-2\pi (n+1)w )})\\
~~~~~=-\sum_{n=0}^{\infty} log(1- e^{-2 \pi \frac{\theta}{w} } e^{\frac{-2\pi (n)}{w} )})+log(1- e^{2 \pi \frac{\theta }{w}} e^{\frac{-2\pi (n+1)}{w} )})\\~~~~~~~~~~~~~~~-\frac{\pi w}{6}+\frac{\pi}{w}(\theta^2-\theta+\frac{1}{6})-\frac{\pi i}{2}+\pi i \theta$
\vspace{0.25cm}
\\
Relocating the terms,
\vspace{0.25cm}
\\
$\displaystyle ~~~~~~~ \sum_{n=1}^{\infty} log(1- e^{\frac{2\pi\theta }{w}} e^{\frac{-2\pi n}{w} })+\sum_{n=1}^{\infty} log(1- e^{ \frac{-2\pi \theta}{w} } e^{\frac{-2\pi (n-1)}{w} })\\~~~~~~~=\sum_{n=1}^{\infty} log(1- e^{2 \pi i \theta } e^{-2\pi nw })+\sum_{n=1}^{\infty} log(1- e^{-2 \pi i \theta } e^{-2\pi (n-1)w })\\~~~~~~~~~~~~~~ -\frac{\pi}{6}\left(w-\frac{1}{w}\right) +\frac{\pi\theta^2}{w} -\frac{\pi \theta}{w}-\frac{\pi i}{2}+\pi i \theta.~~~~~~~~~~~~~~~~~~~~~~~~~~~~~~~~(20)$
\vspace{0.25cm}
\\
This is exactly (19) if we let $\theta=z$ and $w=v$. This proves the transformation law when $k=1$.
\vspace{0.3cm}
\\For $k>1$ we let 
\[\alpha=\frac{\mu}{k}~~where~~1\leq\mu\leq k-1,\]
and writing $h\mu=qk+\phi$ we have again
\[\beta=\frac{\phi}{k}~~where~~1\leq \phi\leq k-1.\]
Note that $\phi\equiv h\mu~(mod~k)$ so $-H\phi\equiv -Hh\mu\equiv \mu~(mod~k)$, and hence $-H\phi/k\equiv \mu/k~(mod~1)$. Therefore \[\alpha= \mu/k\equiv -H\phi/k~(mod~1)\] \[\beta= \phi/k\equiv h\mu/k~(mod~1).\]
Plugging in (18) where again $w=v$ and $\theta=z$ we obtain\\
$\displaystyle ~~~~~~ \sum_{n=0}^{\infty}log(1-e^{\frac{-2\pi z}{v}}e^{-2\pi((n+\beta)v^{-1}+i\alpha)})+\sum_{n=0}^{\infty} log(1-e^{\frac{2\pi z}{v}}e^{-2\pi((n+1-\beta)v^{-1}-i\alpha)})\\
~~~~~~=\sum_{n=0}^{\infty} log(1-e^{2\pi i z}e^{-2\pi((n+\alpha)v-i\beta)})+\sum_{n=0}^{\infty} log(1-e^{-2\pi i z}e^{-2\pi((n+1-\alpha)v+i\beta)})\\~~~~~~~~-\pi v\left(\alpha^2-\alpha+\frac{1}{6}\right)+\frac{\pi}{v}\left(\left(\beta+z\right)^2-(\beta+z)+\frac{1}{6}\right)+2\pi i \left(\alpha-\frac{1}{2}\right)\left(\beta-\frac{1}{2}\right)\\~~~~~~+2\pi zi \left(\alpha-\frac{1}{2}\right).$
\vspace{0.25cm}
\\
Using $\alpha\equiv -H\phi/k~(mod~1)$ and $\beta\equiv h\mu/k~(mod~1)$, we obtain
\vspace{0.25cm}
\\
$\displaystyle ~~~~~~~~~\sum_{n=0}^{\infty}log(1-e^{\frac{-2\pi z}{v}}e^{-2\pi((n+\phi/k)v^{-1}-i\frac{H\phi}{k})})+\sum_{n=0}^{\infty} log(1-e^{\frac{2\pi z}{v}}e^{-2\pi((n+1-\phi/k)v^{-1}+i\frac{H\phi}{k})})\\
~~~~~~~~=\sum_{n=0}^{\infty} log(1-e^{2\pi i z}e^{-2\pi((n+\mu/k)v-i\frac{h\mu}{k})})+\sum_{n=0}^{\infty} log(1-e^{-2\pi i z}e^{-2\pi((n+1-\mu/k)v+i\frac{h\mu}{k})})\\~~~~~~-\pi v\left(\left(\frac{\mu}{k}\right)^2-\frac{\mu}{k}+\frac{1}{6}\right)+\frac{\pi}{v}\left(\left(\frac{\phi}{k}+z\right)^2-(\frac{\phi}{k}+z)+\frac{1}{6}\right)\\~~~~~~~~~+2\pi i \left(\frac{\mu}{k}-\frac{1}{2}\right)\left(\frac{\phi}{k}-\frac{1}{2}\right)+2\pi zi \left(\frac{\mu}{k}-\frac{1}{2}\right).~~~~~~~~~~~~~~~~~~~~~~~~~~~~~~~~~~~~~~~~~~~~~~~(21)$
\vspace{0.27cm}
\\
Note that $log(1-e^{-2\pi(x+mi)})=log(1-e^{-2\pi x})$, i.e it's periodic of period $i$ so the above can be written as\\
$\displaystyle ~~~~~~~\sum_{n=0}^{\infty}log(1-e^{\frac{-2\pi z}{v}}e^{-2\pi(\frac{(nk+\phi)(\frac{1}{v}-iH)}{k})})+\sum_{n=0}^{\infty} log(1-e^{\frac{2\pi z}{v}}e^{-2\pi(\frac{(nk+k+\phi)(\frac{1}{v}-iH)}{k})})\\
~~~~~~~=\sum_{n=0}^{\infty} log(1-e^{2\pi i z}e^{-2\pi\frac{(nk+\mu)(v-ih)}{k})})+\sum_{n=0}^{\infty} log(1-e^{-2\pi i z}e^{-2\pi\frac{(nk+k+\mu)(v-ih)}{k})})\\~~~~~~-\pi v\left(\left(\frac{\mu}{k}\right)^2-\frac{\mu}{k}+\frac{1}{6}\right)+\frac{\pi}{v}\left(\left(\frac{\phi}{k}\right)^2+2z(\frac{\phi}{k})+z^2-(\frac{\phi}{k}+z)+\frac{1}{6}\right)\\~~~~~~+2\pi i \left(\frac{\mu}{k}-\frac{1}{2}\right)\left(\frac{\phi}{k}-\frac{1}{2}\right)+2\pi zi \left(\frac{\mu}{k}-\frac{1}{2}\right)$
\vspace{0.3cm}
\\
Now sum both sides on $\mu$ from $\mu=1,2...k-1$ and also notice that 
\[\left\lbrace nk+\mu, n=0,1,2...; \mu=1,2,...k-1\right\rbrace=\left\lbrace r:r\not \equiv 0~(mod~k)\right\rbrace,\]
and the same goes for the set of number $nk+k-\mu$, and since $\phi\equiv h\mu~(mod~k)$ as $\mu$ runs over the number $1,2,...k-1$ so does $\phi$ but in some other order. Hence we get
\vspace{0.2cm}
\\
$\displaystyle ~~~~~~~\sum_{\substack{r=1 \\ r\not \equiv 0(mod~k)}}^{\infty}log(1-e^{\frac{-2\pi z}{v}}e^{-2\pi r(\frac{(\frac{1}{v}-iH)}{k})})+\sum_{\substack{r=1 \\ r\not \equiv 0(mod~k)}}^{\infty} log(1-e^{\frac{2\pi z}{v}}e^{-2\pi r(\frac{(\frac{1}{v}-iH)}{k})})\\
~~~~~~~=\sum_{\substack{r=1 \\ r\not \equiv 0(mod~k)}}^{\infty} log(1-e^{2\pi i z}e^{-2\pi r\frac{(v-ih)}{k})})+\sum_{\substack{r=1 \\ r\not \equiv 0(mod~k)}}^{\infty} log(1-e^{-2\pi i z}e^{-2\pi r\frac{(v-hi)}{k})})\\~~~~~~~-\pi v\sum_{\mu=1}^{k-1}\left(\left(\frac{\mu}{k}\right)^2-\frac{\mu}{k}+\frac{1}{6}\right)+\frac{\pi}{v}\sum_{\mu=1}^{\infty}\left(\left(\frac{\phi}{k}\right)^2+2z(\frac{\phi}{k})+z^2-(\frac{\phi}{k}+z)+\frac{1}{6}\right) \\~~~~~~~~~+2\pi i \sum_{\mu=1}^{\infty}\left(\frac{\mu}{k}\right)\left(\frac{\phi}{k}-\frac{1}{2}\right)- \pi i \sum_{\mu=1}^{\infty}\left(\frac{\phi}{k}-\frac{1}{2}\right)+2\pi zi \sum_{\mu=1}^{\infty}\left(\frac{\mu}{k}-\frac{1}{2}\right).~~~~~~~~~~~~~~~~(22)$
\vspace{0.25cm}
\\
Checking \cite{apostol}-III, one can see that
\[ \sum_{\mu=1}^{\infty}\left(\frac{\mu}{k}\right)\left(\frac{\phi}{k}-\frac{1}{2}\right)=s(h,k),\]
so ($22$) transforms into\\
$\displaystyle ~~~~~~ \sum_{\substack{r=1 \\ r\not \equiv 0(mod~k)}}^{\infty}log(1-e^{\frac{-2\pi z}{v}}e^{2\pi ir(\frac{(\frac{i}{v}+H)}{k})})+\sum_{\substack{r=1 \\ r\not \equiv 0(mod~k)}}^{\infty} log(1-e^{\frac{2\pi z}{v}}e^{2\pi ir(\frac{\frac{i}{v}+H}{k})})\\
~~~~~~~=\sum_{\substack{r=1 \\ r\not \equiv 0(mod~k)}}^{\infty} log(1-e^{2\pi i z}e^{2\pi i r\frac{(vi+h)}{k})})+\sum_{\substack{r=1 \\ r\not \equiv 0(mod~k)}}^{\infty} log(1-e^{-2\pi i z}e^{2\pi i r\frac{(vi+h)}{k})})\\~~~~~~~-\pi v\sum_{\mu=1}^{k-1}\left(\left(\frac{\mu}{k}\right)^2-\frac{\mu}{k}+\frac{1}{6}\right)+\frac{\pi}{v}\sum_{\mu=1}^{k-1}\left(\left(\frac{\phi}{k}\right)^2+2z(\frac{\phi}{k})+z^2-(\frac{\phi}{k}+z)+\frac{1}{6}\right)\\~~~~~~~+2\pi i s(h,k)-\pi i \sum_{\mu=1}^{k-1}\left(\frac{\phi}{k}-\frac{1}{2}\right)+2\pi zi \sum_{\mu=1}^{k-1}\left(\frac{\mu}{k}-\frac{1}{2}\right).~~~~~~~~~~~~~~~~~~~~~~~~(23)$
\vspace{0.25cm}
\\
Notice that the four sums resemble the desired form, so we intend to look at the residues $g_0$, where\\
$\displaystyle ~~~~~~~ -\pi v\sum_{\mu=1}^{k-1}\left(\left(\frac{\mu}{k}\right)^2-\frac{\mu}{k}+\frac{1}{6}\right)+\frac{\pi}{v}\sum_{\mu=1}^{k-1}\left(\left(\frac{\phi}{k}\right)^2+2z(\frac{\phi}{k})+z^2-(\frac{\phi}{k}+z)+\frac{1}{6}\right)\\~~~~~~~~~+2\pi i s(h,k)-\pi i \sum_{\mu=1}^{k-1}\left(\frac{\phi}{k}-\frac{1}{2}\right)+2\pi zi \sum_{\mu=1}^{k-1}\left(\frac{\mu}{k}-\frac{1}{2}\right)\\~~~~~~~~=-\pi v\left(\frac{1}{k^2}.\frac{k(k-1)(2k-1)}{6}\right)+\pi v\left(\frac{1}{k}.\frac{k(k-1)}{2}\right)-\pi v\left(\frac{k-1}{6}\right)\\~~~~~~~~~~~~~+~2\pi i s(h,k)+\frac{\pi}{v}\left(\frac{1}{k^2}.\frac{k(k-1)(2k-1)}{6}+(2z-1)\left(\frac{k-1}{2}\right)+(z^2-z+\frac{1}{6})(k-1)\right).~~~~~~~~~~~~~(24)$
\vspace{0.3cm}
\\
Note that 
\[-\pi i \sum_{\mu=1}^{k-1}\left(\frac{\phi}{k}-\frac{1}{2}\right)+2\pi zi \sum_{\mu=1}^{k-1}\left(\frac{\mu}{k}-\frac{1}{2}\right)=0\]
since \[\displaystyle \sum_{\mu=1}^{k-1}\left(\frac{\mu}{k}-\frac{1}{2}\right)=\frac{1}{k}\left(\frac{k(k-1)}{2}\right)-\frac{(k-1)}{2}=0,\] which is the same for
$\displaystyle \sum_{\mu=1}^{k-1}\left(\frac{\phi}{k}-\frac{1}{2}\right)$.\\
After some simplification of the terms, (24) becomes
\vspace{0.25cm}
\\
$\displaystyle ~~~~~~~-\frac{\pi}{6k}(k-1)(2k-1)\left(v-\frac{1}{v}\right)+\frac{\pi}{3}(k-1)\left(v-\frac{1}{v}\right)+\frac{k\pi z^2}{v}-\frac{\pi z^2}{v}\\
~~~~~~~=\frac{\pi}{6}\left(v-\frac{1}{v}\right)-\frac{\pi}{6k}\left(v-\frac{1}{v}\right)+\frac{k\pi z^2}{v}-\frac{\pi z^2}{v}.$
\vspace{0.25cm}
\\
Hence (23) is equivalent to \\
$\displaystyle ~~~~~~~\sum_{\substack{r=1 \\ r\not \equiv 0(mod~k)}}^{\infty}log(1-e^{\frac{-2\pi z}{v}}e^{2\pi ir(\frac{(\frac{i}{v}+H)}{k})})+\sum_{\substack{r=1 \\ r\not \equiv 0(mod~k)}}^{\infty} log(1-e^{\frac{2\pi z}{v}}e^{2\pi ir(\frac{\frac{i}{v}+H}{k})})\\
~~~~~~~~~=\sum_{\substack{r=1 \\ r\not \equiv 0(mod~k)}}^{\infty} log(1-e^{2\pi i z}e^{2\pi i r\frac{(vi+h)}{k})})+\sum_{\substack{r=1 \\ r\not \equiv 0(mod~k)}}^{\infty} log(1-e^{-2\pi i z}e^{2\pi i r\frac{(vi+h)}{k})})\\~~~~~~~+2\pi i s(h,k)+\frac{\pi}{6}\left(v-\frac{1}{v}\right)-\frac{\pi}{6k}\left(v-\frac{1}{v}\right)+\frac{k\pi z^2}{v}-\frac{\pi z^2}{v}.$ ~~~~~~~~~~(25) 
\vspace{0.3cm}
\\
Adding equation (25) to equation (20), which corresponds to the case when $k=1$,
\vspace{0.3cm}
\\
$\displaystyle ~~~~~~~~\sum_{n=1}^{\infty} log(1- e^{\frac{2\pi z }{v}} e^{\frac{-2\pi n}{v} })+\sum_{n=1}^{\infty} log(1- e^{ \frac{-2\pi z}{v} } e^{\frac{-2\pi (n-1)}{v} }) \\ ~~~~~~~=\sum_{n=1}^{\infty} log(1- e^{-2 \pi i z } e^{-2\pi (n-1)v })+\sum_{n=1}^{\infty} log(1- e^{2 \pi i z } e^{-2\pi nv }) -\frac{\pi}{6}\left(v-\frac{1}{v}\right) \\ ~~~~~~~~ +\frac{\pi z^2}{v} -\frac{\pi z}{v}-\frac{\pi i}{2}+\pi i z.$ 
\vspace{0.3cm}
\\
This accounts for the missing $r$ where $r\equiv 0~(mod~k)$ if we write $r=mk$,
then the functional equation becomes\\
$\displaystyle ~~~~~~~~ \sum_{\substack{r=1 \\ r\not \equiv 0(mod~k)}}^{\infty}log(1-e^{\frac{-2\pi z}{v}}e^{2\pi ir(\frac{(H+\frac{i}{v})}{k})})+\sum_{\substack{r=1 \\ r\not \equiv 0(mod~k)}}^{\infty} log(1-e^{\frac{2\pi z}{v}}e^{2\pi i(r-1)(\frac{H+\frac{i}{v}}{k})})\\
~~~~~~~~=\sum_{\substack{r=1 \\ r\not \equiv 0(mod~k)}}^{\infty} log(1-e^{2\pi i z}e^{2\pi i r\frac{(h+iv)}{k})})+\sum_{\substack{r=1 \\ r\not \equiv 0(mod~k)}}^{\infty} log(1-e^{-2\pi i z}e^{2\pi i (r-1)\frac{(h+iv)}{k})})\\ ~~~~~~~~~+2\pi i s(h,k)-\frac{\pi}{6k}\left(v-\frac{1}{v}\right)+\frac{\pi k z^2}{v}-\frac{\pi z}{v}+\pi i z -\frac{\pi i }{2}.$
\vspace{0.25cm}
\\This is exactly (17), and this completes the proof of theorem 2 for all $ z \in \mathbf{C}$ by analytic continuation, except at the endpoints $0$ and $\displaystyle \frac{1}{k}$ where we treat them separately. The case for $z=0$ can be trivially tackled and that's why we only treat $\displaystyle z=\frac{1}{k}=\frac{1}{c}$, i.e 
\[~~~~~~~~~~\theta_1\left(\frac{1/c}{c\tau+d},\frac{a\tau+b}{c\tau+d}\right)=\epsilon_1(A)\left(-i(c\tau+d)\right)^{1/2}e^{\frac{\pi i }{c(c\tau+d)}}\theta_1(1/c,\tau).~~~~~~~~~~~~~~~~~~~~~~~~~~~~~~~~~~~~~~~~~~~(26)\]
In order to prove that we do a slight trick. Note that for $z=2(c\tau+d)+1/c$ the theorem holds since $c\tau+d\neq 0$. Hence,
\[~~~~~~~\theta_1\left(\frac{2(c\tau+d)+1/c}{c\tau+d},\frac{a\tau+b}{c\tau+d}\right)=\epsilon_1(A)\left(-i(c\tau+d)\right)^{1/2}e^{\frac{\pi ic(2(c\tau +d)+1/c)^2 }{(c\tau+d)}}\theta_1(2(c\tau+d)+1/c,\tau).\]
i.e,
\[~~~~~~~\theta_1\left(2+\frac{1/c}{c\tau+d},\frac{a\tau+b}{c\tau+d}\right)=\epsilon_1(A)\left(-i(c\tau+d)\right)^{1/2}e^{\frac{\pi ic(2(c\tau +d)+1/c)^2 }{c\tau+d}}\theta_1(2(c\tau+d)+1/c,\tau).\]
Note now that $\theta_1$ is periodic of period 2, hence
\[~~~~~~~\theta_1\left(\frac{1/c}{c\tau+d},\frac{a\tau+b}{c\tau+d}\right)=\epsilon_1(A)\left(-i(c\tau+d)\right)^{1/2}e^{\frac{\pi ic(2(c\tau +d)+1/c)^2 }{c\tau+d}}\theta_1(2c\tau+\frac{1}{c},\tau).~~~~~~~~~~~~~~~~~~~~~~~(27)\]
We need to show that (27) is equivalent to (26), i.e 
\[~~~~~~~~~~~~~~~~~~~~~~~~~~~~~~~~~~~~~~~e^{\frac{\pi i }{c(c\tau+d)}}\theta_1\left(\frac{1}{c},\tau\right)=e^{\frac{\pi ic(2(c\tau +d)+1/c)^2 }{c\tau+d}}\theta\left(2c\tau+\frac{1}{c},\tau\right).~~~~~~~~~~~~~~~~~~~~~~~~~~~~~~~~~~~(28)\]
Note that 
\[e^{\frac{\pi i(2(c\tau +d)+1/c)^2 }{c\tau+d}}=e^{4\pi i c(c\tau+d)}.e^{\frac{\pi i }{c(c\tau+d)}}=e^{4\pi i c^2\tau}.e^{\frac{\pi i }{c(c\tau+d)}}.\]
Taking $\displaystyle e^{\frac{\pi i }{c(c\tau+d)}}$ from both sides in (28), we just have to prove 
\[~~~~~~~~~~~~~~~~~~~~~~~~~~~~~~~~~~~~~~~~~\theta_1\left(\frac{1}{c},\tau\right)=e^{4\pi i c^2\tau}\theta_1\left(2c\tau+\frac{1}{c},\tau\right).~~~~~~~~~~~~~~~~~~~~~~~~~~~~~~~~~~~~~~~~~~~~~~~~~~~~~~~~~~~~~~~(29)\]
Using the relations found in (\cite{theta vocabulary}, 5), given by
\[\theta_1(u+\tau,\tau)=-e^{-\pi i(2u+\tau)}\theta_1(u,\tau).\]
Applying this again, \\ \\
\(~~~~~~~\displaystyle \theta_1(u+2\tau,\tau)=\theta_1(u+\tau+\tau,\tau)=e^{-\pi i (2(u+\tau)+\tau)}\theta_1(u+\tau,\tau)  \\ \displaystyle ~~~~~~~~~~~~~~~~~~~~~~~~~~~~~~~~~~~~~~~~~~~~~~~~=e^{-\pi i(4u+4\tau)}\theta_1(u,\tau)\) \\
Proceeding in this fashion iteratively we obtain
\[\theta_1(u+m\tau,\tau)=(-1)^me^{-\pi i(2mu+m^2\tau)}\theta_1(u,\tau).\]
Setting $\displaystyle u=\frac{1}{c}$ and $m=2c$, we obtain 
\[\theta_1\left(2c\tau+\frac{1}{c},\tau\right)=e^{-4\pi i c^2\tau}\theta_1\left(\frac{1}{c},\tau\right).\]
Plugging in (29), this completes the proof for $z=1/k$ and thus proves theorem 2.
\begin{center}
ACKNOWLEDGMENT  
\end{center} 
We would like to express our deepest gratitude to our advisor Professor Wissam Raji. We also like to thank the Department of Mathematics and the Center of Advanced Mathematical Sciences
(CAMS) at the American university of Beirut (AUB) for the guidance and support we are receiving from the summer research camp (SRC).

Department of Mathematics, American University of Beirut, Beirut, Lebanon
\\
\textit{E-mail address: mmm133@mail.aub.edu}\\
\textit{E-mail address: ays11@mail.aub.edu}
\end{document}